\documentclass[10pt]{article}
\usepackage{amssymb}
\usepackage{amsfonts}
\usepackage{graphicx}
\usepackage{amsmath}
\usepackage{hyperref}
\usepackage{amssymb}
\usepackage{float,color,ulem}
\usepackage{epsf,epsfig,subfigure}
\usepackage{float,color,ulem}
\usepackage[T1]{fontenc}
\usepackage[latin1]{inputenc}
\usepackage{amsmath, amssymb}
\usepackage{listings}
\usepackage{mathrsfs}

\newtheorem{theorem}{Theorem}[section]

\newtheorem{assumption}[theorem]{Assumption}

\newtheorem{corollary}[theorem]{Corollary}

\newtheorem{definition}[theorem]{Definition}
\newtheorem{example}[theorem]{Example}

\newtheorem{lemma}[theorem]{Lemma}

\newtheorem{remark}[theorem]{Remark}

\newenvironment{proof}[1][Proof]{\noindent\textit{#1.} }{\hfill \rule{0.5em}{0.5em}}

\renewcommand{\d}{{\rm d}}

\newcommand{\R}{{\mathbb R}}
\newcommand{\N}{{\mathbb N}}                  
\newcommand{\F}{{\mathbb F}}

\def\Om{\Omega}

\def\t{\tau}

\def\calL{{\mathcal{L}}}

\def\calC{{\mathcal{C}}}

\def\calF{{\mathcal{F}}}

\def\R{\mathbb R}
\def\calN{\mathcal{N}}
\def\P{\mathbb P}

\def\N{\mathbb N}

\def\E{\mathbb E}

\def\F{\mathbb F}

\begin{document}

\title{\textbf{Abstract Parabolic Equations with boundary white noise: an integrated semigroup approach}}
\author{\textsc{Arnaud Ducrot*, Fatima Zahra Lahbiri**} \\
{\small \textit{*Normandie Univ., UNIHAVRE, LMAH, FR-CNRS-3335, ISCN, 76600 Le Havre, France}}\\
{\small	\textit{**FernUniversit\"{a}t in Hagen, Faculty of Mathematics and Computer Science,}} \\
{\small \textit{Chair of Applied Stochastics, D-58084 Hagen, Germany}}\\
		{\small \textit{email: arnaud.ducrot@univ-lehavre.fr, fatima-zahra.lahbiri@fernuni-hagen.de}}
		}
		
\maketitle

\begin{abstract}
In this paper, we study the existence of solution for stochastic evolution equations with almost sectorial operators and possibly a non dense domain. Such problems cover several types of evolution equations, we are interested here in particular in evolution equations with non-homogenous boundary conditions of white noise type. We obtain the existence and uniqueness of mild solutions in state space using the integrated semigroup theory. The results are applied to stochastic parabolic equations with Neumann boundary conditions.

\vspace{0.2in}\noindent \textbf{Key words} Integrated semigroup, Stochastic evolution equation, Parabolic equations, White-noise.

\vspace{0.1in}\noindent \textbf{MSC:}  {\bf 47D06, 37H05, 47D62, 60H40}.
\end{abstract}

\section{Introduction}\label{sec:00}
Given $H$ a separable Hilbert space and closely defined linear operator $A:D(A)\subset H\to H$ with possibly $\overline{D(A)}\neq H$, we consider the stochastic Cauchy problem
\begin{align}\label{Cauchy-pb}
\begin{cases}
dX(t)=AX(t)dt+dW(t), \quad t\in[0,\tau],\cr X(0)=\xi\in\overline{D(A)}.
\end{cases}
\end{align}
Here $\dot{W}(t)$, $t\geq0$ represents a white noise process defined on a probability space $(\Omega,\calF,\P)$, with natural filtration $\F=(\calF_t)_{t\geq0}$.  $\dot{W}$ is a linear transformation from $\mathcal{H}:=L^{2}(0,\tau;H)$ to $L^2(\Omega,\calF,\P)$ with values being Gaussian, zero mean, random variables. The cylindrical Wiener process and white-noise process are defined respectively on $\mathcal{H}$ if
\begin{align*}
  \mathbb{E}(W(\varphi)W(\psi)&=\int_{0}^{\tau}\int_{0}^{\tau}t\wedge s\langle\varphi(t),\psi(s)\rangle dtds \\
   \mathbb{E}(\dot{W}(\varphi)\dot{W}(\psi)&=\int_{0}^{\tau}\langle\varphi(s),\psi(s)\rangle ds.
\end{align*}
If $\{\eta_k\}_{k=1}^{\infty}$ be a sequence of independent normalized Gaussian random variables and $\{\varphi_{k}\}_{k=1}^{\infty}$ be a complete orthonormal basis in $\mathcal{H}$. Then the cylindrical Wiener process $W$ and the white-noise $\dot{W}$ can be represented as a $\P$-a.s. convergent series
\begin{align*}
W(\varphi)&=\sum_{k=1}^{\infty}\eta_k\int_{0}^{\tau}\left\langle\int_{0}^{t}\varphi_k (s)ds,\varphi(t)\right\rangle dt,\quad\varphi\in\mathcal{H},\\
\dot{W}(\varphi)&=\sum_{k=1}^{\infty}\eta_k\int_{0}^{\tau}\langle\varphi_k (t)ds,\varphi(t)\rangle dt,\quad\varphi\in\mathcal{H}.
\end{align*}
If the space $\mathcal{H}$ is embedded in a Hilbert $\tilde{\mathcal{H}}$ with Hilbert-Schmidt embedding then the random series
\begin{align*}
W(t)=\sum_{k=1}^{\infty}\eta_k\int_{0}^{t}\varphi_k (s)ds\quad\dot{W}(t)=\sum_{k=1}^{\infty}\eta_k\varphi_k (t), \quad t\in [0,\tau],
\end{align*}
are almost surely convergent in $\tilde{\mathcal{H}}$. Fore more information on the above concepts as well as their motivation, we refer to \cite{Balakrishnan},\cite{Daprato-93},\cite{Da-Za}.\\
The first objective of this paper is to investigate the existence of the integrated solutions of the stochastic evolution equation \eqref{Cauchy-pb} when $\overline{D(A)}\neq H$, and the linear operator $A$ is almost sectorial. As far as we know, there is a huge literature related to the case where $A$ is a densely defined Hille-Yosida operator (see for instance \cite{Da-Za,Da-Za-96,Duncan-09,Duncan-12} and the references cited therein), here we show that this is not necessary (in a certain sense) to solve problem \eqref{Cauchy-pb}, such results are even more general than those known in the case $\overline{D(A)}=H$.\\
In the absence of the noise (i.e. in the deterministic case) operators with non-dense domains frequently arise in several problems. Some examples are given by adjoint operators of generators on nonreflexive Banach spaces (see \cite{Daprato-67}). Such operators also occur in various biological models such as structured population models (see \cite{Thieme-90}). To guide our investigation, we consider the non-homogenous Cauchy problem with a non-dense domain reads for $f\in L^{1}(0,\tau;H)$ as
\begin{equation}\label{Cauchy}
\begin{cases}
\frac{du(t)}{dt}=Au(t)+f(t),\quad t\in (0,\tau],\\
u(0)\in\overline{D(A)}.
\end{cases}
\end{equation}
The problem \eqref{Cauchy} has been firstly investigated by Da Prato and Sinestrari \cite{Da-Sin} when the linear operator $A$ satisfies the Hille-Yosida property, while integrated semigroup concept was also introduced by Arendt \cite{Arendt-87,Arendt-Is} and Arendt et al \cite{A-B-H} to solve evolution equations with non-densely defined Hille-Yosida operators.
In \cite{Frank-88}, Neubrander kept using of integrated semigroups to study the well-posdness of damped second order Cauchy problems as long as the resolvent of $A$ is polynomially bounded on a certain region. This theory of integrated semigroup was further developed by Magal and Ruan (see \cite{Magal-Ruan,Magal-Ruan07} and the references therein) to study linear and semi-linear abstract Cauchy problems in which the closed linear operator is non-densely defined and satisfies a weak Hille-Yosida property.
This theory of integrated semigroup has been successfully applied to deal with delay differential equations and age structured problem in $L^p$ spaces with $p\in (1,\infty)$ and operators are no longer of Hille-Yosida type.
Let us also notice that this theory has also be used to study parabolic equations with non-homogeneous boundary conditions and nonlinear boundary conditions. We refer the reader to \cite{Ducrot-Magal-Prevost} for results on the semigroups generated by almost sectorial operators (see Assumption \ref{Assumption1.1} below) and to \cite{Ducrot-Magal-pise} for applications to semi-linear problems and parabolic equation with nonlinear and nonlocal boundary conditions as a special case.

In the stochastic setting, Neamtu \cite{Neamtu-20} used the theory of integrated semigroups to study the well-posedness of Stratonovich evolution equations and to construct random stable/unstable manifolds using the Lyapunov-Perron method. Recently, Li and Zeng \cite{Li-Zeng 22} developed a similar idea to construct a center manifold for the ill-posed Stratonovich stochastic evolution equations with a non-dense domain. As far as we know there are few results regarding stochastic evolution equations of It\^{o} type, which is the main focus of this work.\\
We aim to investigate the existence of integrated solutions of \eqref{Cauchy-pb} when $A$ is almost sectorial operator using integrated semigroups. When dealing with parabolic equations (densely defined or not), it is usually assumed that the operator $A$ is a sectorial elliptic operator. This operator property usually holds true when considering elliptic operators in Lebesgue spaces or H\"{o}lder spaces and together with homogeneous boundary conditions. As pointed out by Lunardi in \cite{Lunardi}, this property does no longer hold true when dealing with such operators in some more regular spaces. Typical examples of non-sectorial but almost sectorial operators may also arise when dealing with parabolic equations with non-homogeneous boundary conditions, such as stochastic parabolic delay differential equations, stochastic parabolic boundary control equations and parabolic evolution equations with non-homogenous boundary conditions of white noise type. In our setting, we are interested specifically in evolution equations with non-homogeneous boundary white noise, see Section \ref{sec:3}. This type of application was first investigated by Da-prato and Zabczyk \cite{Daprato-93,Da-Za-96} who used an extrapolation semigroup approach to handle parabolic problems with both Dirichlet and Neumann boundary conditions. They obtained the existence and uniqueness of a mild solution in a suitable extrapolation space larger than the natural state space. Maslowski \cite{Maslowski-95} used a similar approach, obtaining an existence and uniqueness result for mild solutions of the problem by a fixed point argument. See also \cite{Alb-Roz93,Deb-all07, Gol-Pes21,Sch-Ver11} and \cite{Sowers} for further contributions to problems with boundary noise. Results regarding stochastic evolution equations of It\^{o}'s type with homogenous boundary conditions can be looked up to the work of Hadd and Lahbiri \cite{Lahbiri22,Lahbiri-population,Lahbiri21}.\\
This note aims to investigate the well-posed of the stochastic abstract Cauchy problem \eqref{Cauchy-pb} when the linear operator $A$ is possibly non-densely defined and satisfies the following assumption.
\begin{assumption}
\label{Assumption1.1} Let $A:D(A)\subset H\rightarrow H$ be a linear operator
on the Hilbert space $\left( H,\left\langle \cdot,\cdot\right \rangle_H\right).$ We assume that

\begin{itemize}
\item[(a)] the operator $A_{0}$, the part of $A$ in $\overline{D(A)}$, is the
infinitesimal generator of an
analytic semigroup of bounded linear operators on $\overline{D(A)}$ that is denoted by $\left\{ T_{A_{0}}(t)\right\} _{t\geq 0}$.

\item[(b)] There exist $\omega \in \mathbb{R}$ and $p^{\ast }\in \left[
1,+\infty \right)$ such that $\left( \omega ,+\infty \right) \subset \rho
\left( A\right) $, the resolvent set of $A$, and
\begin{equation}
\underset{\lambda \rightarrow +\infty }{\limsup }\;\lambda ^{\frac{1}{%
p^{\ast }}}\left\Vert \left( \lambda I-A\right) ^{-1}\right\Vert _{\mathcal{L%
}\left( H\right) }<+\infty.  \label{1.2}
\end{equation}
\end{itemize}
\end{assumption}
Using the above assumption and the results in \cite{Ducrot-Magal-Prevost}, the linear operator $A$ becomes the generator of an analytic integrated semigroup $\{S_A(t)\}_{t\geq 0}$ on $H$.
Thus this family of linear operators turns out to be differentiable with respect to $t>0$ and the derivative $\left\{\frac{dS_A(t)}{dt}\right\}_{t>0}$ forms a semigroup on $H$ which possibly not strongly continuous at $t=0$ due to a singularity characterized by $p^*$ at this point $t=0$.
Roughly speaking we prove that $S_A$ belongs for some $\tau>0$ to $W^{1,p}(0,\tau;\mathcal{H})$ with $p$ sufficiently large then the stochastic convolution integral
\begin{equation*}
\int_0^t \frac{dS_A(t-s)}{ds}dW(s)
\end{equation*}
is well defined almost surely in the state space $H$ and plays the role of the constant variation formula to define mild solutions of \eqref{Cauchy-pb}, as the standard variation of constants formula is not applicable for the case in point (see Theorem \ref{main-result2}).
The result is then applied to parabolic equations with boundary noise, illustrated by the study of a heat equation with a non-homogeneous Neumann boundary condition involving a white noise.\\

Using the results proved by Ducrot et al. in \cite{Ducrot-Magal-Prevost} (see Proposition 3.3 in that paper), this above set of assumptions can be reformulated using the notion of almost sectorial operators and, this re-writes as follows.\\
The Assumption \ref{Assumption1.1} is
satisfied if and only if the two following conditions are satisfied:
\begin{itemize}
\item[(a)] $A_{0}$ is a sectorial operator.
\item[(b)] $A$ is a $\frac{1}{p^{\ast }}-$almost sectorial operator.
\end{itemize}
The definition of $\alpha-$almost sectorial operator is recalled in Definition \ref{DE1.2} Section \ref{sec:1}.\\
Some classical notations on stochastic analysis is also required. Let $\mathscr{X}$ is be Banach space and a sub-$\sigma$- algebra $\calN$ of $\calF$, denote $L^{2}_{\calN}(\Omega;H)$ the set of all $\calN$-measurable ($\mathscr{X}$-valued) random variables $\zeta:\Omega\to\mathscr{X}$ with $\E|\zeta|^{2}_{\mathscr{X}}<\infty$. Next, if $\t\in(0,\infty]$ is a real number, we denote
\begin{align*}
L^2_{\F}(0,\t;\mathscr{X}) :=\left\{\zeta:[0,\t]\times \Om\to \mathscr{X}:\zeta \;\text{is}
\;\F-\text{adapted and}\; \int^\t_0 \mathbb{E}\|\zeta(t)\|_{\mathscr{X}}^2<\infty\right\}
\end{align*}
and the space of  mild solutions for stochastic equations
\begin{align*}
& \mathcal{C}_{\F}\left(0,\t;L^2(\Om,\mathscr{X})\right)\cr & \quad :=\left\{\zeta:[0,\t]\times \Om\to \mathscr{X}:\zeta \;\text{is}\;\F-\text{adapted and}\; t\mapsto\left(\mathbb{E}\|\zeta(t)\|_{\mathscr{X}}^2\right)^{\frac{1}{2}}\;\text{is continuous}\right\}.
\end{align*}
Next, let $\mathfrak{X}_{1}$ and $\mathfrak{X}_{2}$ be two separable Hilbert spaces with an orthonormal basis $\{e_{k}\}_{k=1}^{\infty}$ in $\mathfrak{X}_{1}$, then the space of Hilbert-Schmidt operators from $\mathfrak{X}_{1}$ to $\mathfrak{X}_{2}$ is defined as
\begin{equation}\label{Hilbert-Schmidt}
\calL_{2}(\mathfrak{X}_{1},\mathfrak{X}_{2})=\left\{T\in\calL(\mathfrak{X}_{1},\mathfrak{X}_{2}):\sum_{k\in\N} \|Te_k\|^{2}_{\mathfrak{X}_{2}} <+\infty\right\}.
\end{equation}
It is well known (see \cite{Da-Za},\cite{Gawa-Mand},\cite{Schatten-13}) that $\calL_{2}(\mathfrak{X}_{1},\mathfrak{X}_{2})$ equipped with the norm
\begin{equation*}
\|T\|_{\calL_{2}(\mathfrak{X}_{1},\mathfrak{X}_{2})}=\left(\sum_{k=1}^{\infty} \|Te_k\|^{2}_{\mathfrak{X}_{2}}\right)^{\frac{1}{2}}
\end{equation*}
is a Hilbert space. Since the Hilbert spaces $\mathfrak{X}_{1}$ and $\mathfrak{X}_{2}$ are separable, the space $\calL_{2}(\mathfrak{X}_{1},\mathfrak{X}_{2})$ is also separable, as Hilbert-Schmidt operators are limits of sequences of finite-dimensional linear operators. Note that the number $\|T\|_{\calL_{2}(\mathfrak{X}_{1},\mathfrak{X}_{2})}$ is independent of the choice of orthonormal basis $\{e_{k}\}_{k=1}^{\infty}$ in $\mathfrak{X}_1$.\\

This paper is organized as follows: In Section \ref{sec:1}, we recall the necessary material that will be needed in this work about analytic integrated semigroups. Section \ref{sec:2} is devoted to the existence and the uniqueness of integrated solutions of \eqref{Cauchy-pb}. Finally Section \ref{sec:3} presents applications to parabolic problems with boundary noise. An illustrative example with a stochastic heat equation with Neumann boundary conditions is given.
\section{Preliminary Material On Analytic Integrated Semigroups}\label{sec:1}
In this section we present some materials on linear equations and recall some important results that will be used in the sequel.
Let $X$ and $Z$ be two Banach spaces. We denote by $\mathcal{L}\left( X,Z\right) $
 the space of bounded linear operators from $X$ into $Z$ and by $%
\mathcal{L}\left( X\right) $ the space $\mathcal{L}\left( X,X\right) .$ Let $%
A:D(A)\subset X\rightarrow X$ be a linear operator.
We set
\begin{equation*}
X_{0}:=\overline{D(A)},
\end{equation*}%
and we denote by $A_{0}$, the part of $A$ in $X_{0},$ the linear operator on $X_{0}$
defined by
\begin{equation*}
A_{0}x=Ax,\;\forall x\in D(A_{0}):=\left\{ y\in D(A):Ay\in X_{0}\right\} .
\end{equation*}
Throughout this section we assume that $A$ satisfies Assumption \ref{Assumption1.1} for some $p^*\in [1,\infty)$ and $\omega\in\R$.
Note that it is
easy to check that for each $\lambda >\omega$ one has
\begin{equation*}
D\left( A_{0}\right) =\left( \lambda I-A\right) ^{-1}X_{0}\text{ and }\left(
\lambda I-A_{0}\right) ^{-1}=\left( \lambda I-A\right) ^{-1}\mid _{X_{0}}.
\end{equation*}
 From here on, we define $q^{\ast }\in \left( 1,+\infty \right]$ by
\begin{equation}
q^{\ast }:=\frac{p^{\ast }}{p^{\ast }-1}\text{ }\Leftrightarrow \frac{1}{%
q^{\ast }}+\frac{1}{p^{\ast }}=1,  \label{2.1}
\end{equation}%
wherein $p^{\ast }\geq 1$ is defined in Assumption \ref{Assumption1.1}.\

In order to prepare our ground, we first recall some results for the non-homogeneous Cauchy problems
\begin{equation}
\frac{\d u(t)}{\d t}=Au(t)+f(t),t\geq 0,\text{ }u(0)=x\in \overline{D(A)}.
 \label{2.2}
\end{equation}%
To that aim let us recall the following definition.
\begin{definition}[Integrated solution]
\label{DE2.1}
Let $f\in L^{1}\left( 0,\tau ;X\right)$ be a given function for
some given $\tau >0$. A map $us\in C\left( \left[ 0,\tau \right] ,X\right) $ is
said to be \textbf{an integrated solution} of the Cauchy problem \eqref{2.2} on $%
[0,\tau ]$ if the two following conditions are satisfied:
\begin{equation*}
\begin{split}
&\displaystyle \int_{0}^{t}u(s)\d s\in D(A),\forall t\in \left[ 0,\tau \right] ,\\
&\text{and }\\
&\displaystyle u(t)=x+A\int_{0}^{t}u(s)\d s+\int_{0}^{t}f(s)\d s,\;\forall t\in [0,\tau].
\end{split}
\end{equation*}
\end{definition}
In order to go further recall that $\omega _{0}(A_{0})$ the growth rate of the semigroup $\left\{
T_{A_{0}}(t)\right\} _{t\geq 0}$ is defined by
\begin{equation*}
\omega _{0}(A_{0}):=\lim_{t\rightarrow +\infty }\frac{\ln \left( \left\Vert
T_{A_{0}}(t)\right\Vert _{\mathcal{L}\left( X_0\right) }\right)
}{t}.
\end{equation*}%
Since $p^{\ast }\neq +\infty$, one has $\left\Vert \left( \lambda I-A\right) ^{-1}\right\Vert _{\mathcal{L}\left(
X\right) }\to 0$ as $\lambda \rightarrow +\infty$ and by using the Lemma 2.1\ in Magal and Ruan \cite%
{Magal-Ruan07},\ we deduce that
\begin{equation*}
\overline{D(A)}=\overline{D(A_{0})}.
\end{equation*}%
Since by assumption $\rho \left( A\right) \neq \emptyset ,$ it is follows
that (see Magal and Ruan \cite[Lemma 2.1]{Magal-Ruan09a})
\begin{equation*}
\rho \left( A\right) =\rho \left( A_{0}\right) .
\end{equation*}%
This in particular yields
\begin{equation*}
\left( \omega _{0}(A_{0}),+\infty \right) \subset \rho \left( A\right).
\end{equation*}
Next the integrated semigroup $\left\{ S_{A}(t)\right\} _{t\geq 0}$ generated by $%
A$ is the family of bounded linear operator on $X$ defined for all $\lambda \in \left( \omega _{0}(A_{0}),+\infty \right) $ by
\begin{equation}
S_{A}(t)=\left( \lambda I-A_{0}\right) \int_{0}^{t}T_{A_{0}}(s)\d s\left(
\lambda I-A\right) ^{-1},  \label{2.3}
\end{equation}%
where $(T_{A_{0}})_{t\geq 0}$ is the $\mathcal{C}_0$-semigroup generated by $A_0$.\

The relationship between the integrated semigroups $\left\{ S_{A}(t)\right\}
_{t\geq 0}$, and the semigroup, used in paticular by Lunardi in \cite{Lunardi}, comes from
the fact that the map $t\rightarrow S_{A}(t)$ is continuously differentiable from $%
\left( 0,+\infty \right) $ into $\mathcal{L}\left( X\right) $, and
that the family
\begin{equation}
T_A(t):=\frac{\d S_{A}(t)}{\d t}=\left( \lambda I-A_{0}\right) T_{A_{0}}(t)\left(
\lambda I-A\right) ^{-1},\text{ for }t>0,\text{ and }T_A(0)=I,
\label{2.4}
\end{equation}%
defined a semigroup of bounded linear operators on $X$. However it has to be noted that when $A$ is not
densely defined then the family $\left\{ T_A(t)\right\} _{t\geq 0}$ of bounded
linear operator on $X$ is not strongly continuous at $t=0$.\\
For completeness, we also recall that the analyticity of $t\rightarrow S_{A}(t)$ and $%
t\rightarrow T_A(t),$ follows from the formula
\begin{equation*}
S_{A}(t)=\left( \mu I-A_{0}\right) \int_{0}^{t}T_{A_{0}}(l)\d l\left( \mu
I-A\right) ^{-1},\text{ and }T_A(t)=\int_{\Gamma }e^{\lambda t}(\lambda
-A)^{-1}\d\lambda ,
\end{equation*}%
where $\mu >\omega _{0}\left( A_{0}\right) ,$ and $\Gamma $ is the path $%
\omega+\left\{ \lambda \in \mathbb{C}:\left\vert \arg \left( \lambda
\right) \right\vert =\eta ,\left\vert \lambda \right\vert \geq r\right\}
\cup $ $\left\{ \lambda \in \mathbb{C}:\left\vert \arg \left( \lambda
\right) \right\vert \leq \eta ,\left\vert \lambda \right\vert =r\right\} $,
oriented counterclockwise for some $r>0$, $\eta \in \left( \frac{\pi}{2},\pi\right) $.\\
Now, let's recall the concept of almost sectorial operators.
\begin{definition}[Almost sectorial operator]
\label{DE1.2} Let $L:D(L)\subset X\rightarrow X$ be a linear operator on the Banach
space $X$ and let $\alpha \in \left( 0,1\right] $ be given. Then $L$ is said to be a $%
\alpha -$\textbf{almost sectorial} operator if there are constants $\widehat{%
\omega }\in \mathbb{R},$ $\theta \in \left( \frac{\pi }{2},\pi \right) ,$
and $\widehat{M}>0$ such that
\begin{itemize}
\item[(i)] $\rho (L)\supset S_{\theta ,\widehat{\omega }}=\left\{ \lambda
\in \mathbb{C}:\lambda \neq \widehat{\omega },\left\vert \arg \left( \lambda-\widehat{\omega }\right) \right\vert <\theta \right\} ,$
\item[(ii)] $\displaystyle \left\Vert \left( \lambda I-L\right) ^{-1}\right\Vert_{\mathcal L(X)} \leq
\frac{\widehat{M}}{|\lambda -\widehat{\omega }|^\alpha},\;\;\forall \lambda \in S_{\theta ,\widehat{\omega }}.$
\end{itemize}
Moreover $L$ is called \textbf{sectorial operator} if $L$ is $1$-almost
sectorial.
\end{definition}

In the context of Assumption \ref{Assumption1.1}, recall also that the fractional powers $\left( \lambda I-A_{0}\right) ^{-\alpha }$ are well defined, for any $\lambda >\omega _{0}(A_{0})$, by
\begin{equation*}
\left( \lambda I-A_{0}\right) ^{-\alpha }=\frac{1}{\Gamma (\alpha )}%
\int_{0}^{+\infty }t^{\alpha -1}T_{\left( A_{0}-\lambda I\right) }(t)\d t,%
\text{ for }\alpha >0,\text{ and }\left( \lambda I-A_{0}\right) ^{0}=I.
\end{equation*}%
Now since $A$ is only assumed to be almost sectorial, the fraction powers of $\left(
\lambda I-A\right) ^{-\alpha }$ are not defined for any $\alpha >0$ but for $\alpha$ large enough. More
precisely, we have following result (see \cite{periago} or \cite[Lemma 3.7]{Ducrot-Magal-Prevost}.

\begin{lemma}
\label{LE2.2}Let Assumption \ref{Assumption1.1} be satisfied. The fractional
power $\left( \lambda I-A\right) ^{-\alpha }\in \mathcal{L}\left( X\right) $
is well defined for each $\alpha \in \left( \dfrac{1}{q^{\ast }},+\infty
\right)$ and $\lambda >\omega _{0}(A_{0})$.\ Moreover one has
\begin{equation*}
\left( \lambda I-A\right) ^{-\alpha }\left( X\right) \subset \overline{D(A)},
\end{equation*}%
and the following properties are satisfied:
\begin{itemize}
\item[(i)] $\left( \mu I-A_{0}\right) ^{-1}\left( \lambda I-A\right)
^{-\alpha }=\left( \lambda I-A_{0}\right) ^{-\alpha }\left( \mu I-A\right)
^{-1},\forall \mu >\omega _{0}(A_{0}).$

\item[(ii)] $\left( \lambda I-A_{0}\right) ^{-\alpha }x=\left( \lambda
I-A\right) ^{-\alpha }x,\;\forall x\in \overline{D(A)}=X_0.$

\item[(iii)] For each $\alpha \geq 0,$ $\beta >\dfrac{1}{q^{\ast }},$%
\begin{equation*}
\left( \lambda I-A_{0}\right) ^{-\alpha }\left( \lambda I-A\right) ^{-\beta
}=\left( \lambda I-A\right) ^{-\left( \alpha +\beta \right) }.
\end{equation*}
\end{itemize}
\end{lemma}

Now observe that since $\left( \lambda I-A\right) ^{-\alpha }$ and $\left( \mu
I-A\right) ^{-1}$ commute, it follows that $\left( \lambda I-A\right) ^{-\alpha }$ commutes with $S_{A}(t)$ and $T_{A_{0}}(t).$ This in particular yields
\begin{equation*}
S_{A}(t)=\left( \lambda I-A_{0}\right) ^{\alpha
}\int_{0}^{t}T_{A_{0}}(s)ds\left( \lambda I-A\right) ^{-\alpha }
\end{equation*}%
for any $\alpha \in \left( \dfrac{1}{q^{\ast }},+\infty \right)$ and for each $%
\lambda >\omega _{0}(A_{0})$.

Let us also observe that for $\alpha \in \left( \dfrac{1}{q^{\ast }},1\right] ,$
\begin{equation*}
\left( \lambda I-A\right) ^{-1}=\left( \lambda I-A_{0}\right) ^{-(1-\alpha
)}\left( \lambda I-A\right) ^{-\alpha }.
\end{equation*}%
Hence, due to \eqref{2.4}, for each $t>0$ we get
\begin{equation*}
\begin{split}
\frac{\d S_{A}(t)}{\d t}&=\left( \lambda I-A_{0}\right) T_{A_{0}}(t)\left( \lambda
I-A\right) ^{-1}\\
&=\left( \lambda I-A_{0}\right) T_{A_{0}}(t)\left( \lambda
I-A_{0}\right) ^{-(1-\alpha )}\left( \lambda I-A\right) ^{-\alpha }
\end{split}
\end{equation*}%
and, since $T_{A_{0}}(t)$ and $\left( \lambda I-A_{0}\right) ^{-(1-\alpha )}$
commute, we also obtain the following expression for the derivative of $S_A$:
\begin{equation}
\frac{\d S_{A}(t)}{\d t}=\left( \lambda I-A_{0}\right) ^{\alpha
}T_{A_{0}}(t)\left( \lambda I-A\right) ^{-\alpha },\forall t>0,\forall
\alpha \in \left( \dfrac{1}{q^{\ast }},1\right] . \label{2.10}
\end{equation}

Now the main tool to deal with integrated solutions for the Cauchy problem relies on the constant variation formula. Hence before coming back to the non-homogeneous problem \eqref{2.2} let us recall the following result.

\begin{theorem}
\label{TH2.3} Let Assumption \ref{Assumption1.1} be satisfied.\ Let $f\in
L^{p}\left( 0,\tau ;X\right) $ with $p>p^{\ast }$. Then the map $%
t\rightarrow \left( S_{A}\ast f\right) (t):=\int_{0}^{t}S_{A}(t-s)f(s)\d s$ is
continuously differentiable, $\left( S_{A}\ast f\right) (t)\in D(A),\;\forall
t\in \left[ 0,\tau \right] ,$ and if we denote by%
\begin{equation}
\left( S_{A}\diamond f\right) (t):=\frac{\d}{\d t}\int_{0}^{t}S_{A}(t-s)f(s)\d s,
 \label{2.5}
\end{equation}%
then%
\begin{equation*}
\left( S_{A}\diamond f\right) (t)=A\int_{0}^{t}\left( S_{A}\diamond f\right)
(s)\d s+\int_{0}^{t}f(s)\d s,\;\forall t\in \left[ 0,\tau \right] .
\end{equation*}%
Moreover for each $\beta \in \left( \dfrac{1}{q^{\ast }},\dfrac{1}{q}\right)
$ (with $\dfrac{1}{q}+\dfrac{1}{p}=1$), each $\lambda >\omega _{0}(A_{0}),$
and each $t\in \left[ 0,\tau \right] ,$ the following holds true
\begin{equation}
\left( S_{A}\diamond f\right) (t)=\int_{0}^{t}\left( \lambda I-A_{0}\right)
^{\beta }T_{A_{0}}(t-s)\left( \lambda I-A\right) ^{-\beta }f(s)\d s,
\label{2.6}
\end{equation}%
and, the following estimate also holds true
\begin{equation}
\left\Vert \left( S_{A}\diamond f\right) (t)\right\Vert \leq M_{\beta
}\left\Vert \left( \lambda I-A\right) ^{-\beta }\right\Vert _{\mathcal{L}%
(X)}\int_{0}^{t}(t-s)^{-\beta }e^{\omega _{A}(t-s)}\left\Vert
f(s)\right\Vert \d s, \label{2.7}
\end{equation}%
wherein $M_{\beta }$ denotes some positive constant, and $\omega _{A}>\omega
_{0}(A_{0})$.
\end{theorem}

\bigskip

By using integrated semigroups, or formula \eqref{2.6}, we derive the
extended variation of constant formula:
\begin{equation}
\left( S_{A}\diamond f\right) (t)=T_{A_{0}}(t-s)\left( S_{A}\diamond
f\right) (s)+\left( S_{A}\diamond f(s+.)\right) (t-s),\forall t\geq s\geq 0.
 \label{2.8}
\end{equation}%
By using the above theorem, and the usual uniqueness result of
Thieme \cite[Theorem 3.7]{Thieme-90a}, one derive the following result.

\begin{corollary}
\label{CO2.4} Let Assumption \ref{Assumption1.1} be satisfied. Let $p\in
\left( p^{\ast },+\infty \right) $ be given. Then for each $f\in
L^{p}\left( 0,\tau ;X\right) $ and for each $x\in X_{0}$ the
Cauchy problem \eqref{2.2} has a unique integrated solution $u\in C\left( %
\left[ 0,\tau \right] ,X_{0}\right) $ that is given by%
\begin{equation}
u(t):=T_{A_{0}}(t)x+\left( S_{A}\diamond f\right) (t),\text{ }\forall t\in %
\left[ 0,\tau \right] . \label{2.9}
\end{equation}
\end{corollary}

\section{Integrated solution}\label{sec:2}
This section aims to present a functional analytic method to establish the mild solutions of the stochastic Cauchy problem \eqref{Cauchy-pb}. We fix $\{\varphi_{j}\}_{j=1}^{\infty}$ an orthonormal and complete basis in $\mathcal{H}$ consisting of smooth functions, and sequences $\{\eta_{j}\}_{j=1}^{\infty}$ of independent $\mathcal{N}(0,1)$ normalized random variables and consider approximation equation
\begin{align}\label{Trans-ap}
\begin{cases}dX_{N}(t)&=\mathscr{A} X_{N}(t)dt+d\mathscr{W}_{N}(t),\quad t\in [0,\tau], \cr X_{N}(0)&=\xi\in \overline{D(A)},\end{cases}
\end{align}
where for any $N\in\N$ and $t\in [0,\tau]$
\begin{equation}\label{white-noise-appro}
  \dot{\mathscr{W}}_{N}(t)=\sum_{j=1}^{N}\eta_{j}\varphi_{j}(t).
\end{equation}
In order to define the stochastic integral, we make the following assumption.
\begin{assumption}
\label{Assumption1.2}
\begin{itemize}
\item[(a)] Assume that
 \begin{equation*}
  \int_{0}^{\tau} \|{S}_{A}(r)\|^{2}_{\calL_{2}(H)}dr=\int_{0}^{\tau}Tr({S}_{A}(r))dr<\infty,
  \end{equation*}
 \item[(b)] As $t\to S_A(t)$ is continuously differentiable, we assume that
 \begin{equation*}
   \int_{0}^{\tau} \left\|\frac{d{S}_{A}(r)}{dr}\right\|^{2}_{\calL_{2}(H)}dr=\int_{0}^{\tau}Tr\left(\frac{d{S}_{A}(r)}{dr}\right)dr<\infty.
\end{equation*}
\end{itemize}
\end{assumption}
Let assumptions \eqref{Assumption1.1} and \eqref{Assumption1.2} hold. Then the following transformation
\begin{equation*}
 ({S}_{A}\ast u)(t)=\int_{0}^{t}{S}_{A}(t-s)u(s)ds,\quad t\in[0,\tau],
\end{equation*}
is linear from $W^{1,2}(0,\tau;{H})$ into $\mathcal{H}$. Then the map $t\to({S}_{A}\ast u)(t)$ is differentiable on $[0,\tau]$ and for each $t\in[0,\tau]$
\begin{align*}
  \frac{d}{dt} ({S}_{A}\ast u)(t)=\int_{0}^{t}\frac{d}{dr}S_{A}(r)u(t-r)dr.
\end{align*}
Next let us show that the process $t\to X_{N}(t)$ is an integrated solution to \eqref{Trans-ap}. Our precise result reads as follows.
\begin{theorem}\label{main-result1}
Let assumptions \eqref{Assumption1.1} and \eqref{Assumption1.2} hold. Then the process $t\mapsto X_N(t)$ satisfies
\begin{align*}
\int_0^t X_N(s,\omega) ds\in D(A),\quad\omega\in\Omega,
\end{align*}
and
\begin{equation*}
X_{N}(t)=\xi+A\int_0^t u(s) ds+\mathscr{W}_N (t),\;\;\forall t\geq 0.
\end{equation*}
\end{theorem}
\begin{proof}
Let us fix $\tau>0$. Then the stochastic Fubini theorem applies and ensures that we have for any $t\in [0,\tau]$
\begin{eqnarray*}
&&\int_0^t \frac{dS_A}{dt}(t-s)\int_0^s d\mathscr{W}_N (u) ds=\int_0^t\int_u^t \frac{dS_A}{dt}(t-s)ds d\mathscr{W}_N (u)\\
&=&\int_0^t S_A(t-u)d\mathscr{W}_N(u)=\int_0^t u(\ell)d\ell.
\end{eqnarray*}
On the other hand we also has for all $t\in [0,\tau]$
\begin{eqnarray*}
&&\int_0^t \frac{dS_A}{dt}(t-s) \int_0^s d\mathscr{W}_N (u)ds=\int_0^t \frac{dS_A}{dt}(t-s)ds  [\mathscr{W}_N (s)-\mathscr{W}_N (0)]\\
&=&\int_0^t \frac{dS_A}{dt}(t-s) \mathscr{W}_N (s)ds=\frac{d}{dt}\left(S_A\ast \mathscr{W}_N \right)(t).
\end{eqnarray*}
Therefore by standard regularity results for the parabolic inhomogeneous Cauchy problem \cite[Theorem 5.3.5]{Lunardi} and by \cite[Theorem 3.11]{Ducrot-Magal-Prevost}, $t\to \frac{d}{dt}\left(S_A\ast \mathscr{W}_N \right)(t,\omega)=\int_0^t \frac{dS_A}{dt}(t-s) \mathscr{W}_N (s,\omega)ds$ belongs to $\mathcal{C}(0,\tau;D(A))$.
\end{proof}
\begin{remark}
An integrated solution always belongs to $H_0:=\overline{D(A)}$.
\end{remark}
Keeping the above remark in mind our next goal is to derive a variation of constant formula for the approximation problem \eqref{Trans-ap}. The difficulty that occurs is that $\{\frac{dS_A(t)}{dt}\}_{t\geq0}$ forms a $\mathcal{C}_0$-semigroup only on $H_0$, however the process $\mathscr{W}_N $ out of $H_0$.
\begin{lemma}\label{lemma-result1}
Let assumptions \eqref{Assumption1.1} and \eqref{Assumption1.2} hold. Then the integrated solution of \eqref{Trans-ap} satisfies the following variation of constant formula
\begin{equation*}
  X_N(t)={T}(t)\xi+\lim_{\lambda\to +\infty} \int_{0}^{t}{T}(t-s)\lambda R(\lambda,{A})d\mathscr{W}_{N}(s).
\end{equation*}
\end{lemma}
\begin{proof}
Let us observe that $X_{N}(t)$ is a mild solution of \eqref{Trans-ap} if and only if $v_{N}(t):=\int_{0}^{t}X_{N}(s) ds$ is a strong solution of
\begin{align}\label{Trans-ap2}
\begin{cases}dv_{N}(t)&={A} v_{N}(t)dt+\xi+\mathscr{W}_{N}(t),\quad t\in [0,\tau], \cr v_{N}(0)&=0.\end{cases}
\end{align}
Note that $\xi\in{H}_{0}$ is a necessary condition for the existence of any kind of solution in our setting. Note that $t\to v_N(t)$ is a strong solution to \eqref{Trans-ap2} and one has
\begin{equation*}
 \int_0^t v_N(s)ds=\left(S_A\ast \mathscr{W}_N \right)(t)=\int_0^t \int_0^s X_N(\ell)d\ell ds\in D(A),
\end{equation*}
Then the solution is given by
\begin{align}\label{approximate-v}
v_{N}(t)=S_{A}(t)\xi+\frac{d}{dt}\left(S_A\ast \mathscr{W}_N \right)(t).
\end{align}
The above Formulation can be done using the same arguments given by \cite{Magal-Ruan}. Now, as $\xi\in{H}_{0}$ and $\dot{\mathscr{W}}_{N}\in W^{1,2}(0,\tau;{H})$, we deduce that the right-hand side is differentiable in $[0,\tau]$ and
\begin{align*}
v^{\prime}_{N}(t)={T}_{0}(t)\xi+\int_{0}^{t}\frac{d}{ds}{S}_{A}(t-s)d\mathscr{W}_{N}(s).
\end{align*}
Multiplying the previous equality by $\lambda R(\lambda,A)$ and letting $\lambda\to\infty$ gives us the desired constant formula due to the fact that $X_N\in\mathcal{C}([0,\tau]; H_0)$ . This can be shown by the same arguments used in the deterministic case, see also \cite[Theorem 1.6]{Magal-Ruan}. The result follows from the fact that $v^{\prime}_{N}(t)=X_{N}(t)$.
\end{proof}
\newline
The following result is the main result of this section.
\begin{theorem}\label{main-result2}
Let assumptions \eqref{Assumption1.1} and \eqref{Assumption1.2} hold and let $X(\cdot)$ be the process given by \eqref{Cauchy-pb}. For any $\xi\in H_{0}$, $X_{N}(\cdot)$ converges to $X(\cdot)$ in $\calC_{\F}([0,\tau],L^{2}(\Omega,H))$ as $N\to +\infty$ and
\begin{align*}
X(t)={T}_{0}(t)\xi+\lim_{\lambda\to +\infty} \int_{0}^{t}{T}_{0}(t-s)\lambda R(\lambda,{A})d{W}(s),
\end{align*}
for any $t\in [0,\tau]$ and $\P$-a.s.
\end{theorem}
\begin{proof}
Let $v_N$ be the process given in \eqref{Trans-ap2}. For $N\in\N$ and $t\in [0,\tau]$, consider $v_{N}(t)$ as in \eqref{approximate-v}. Let $M,N\in\N$ with $M>N$, using H\"{o}lder's inequality we get for any $\{h_{k}\}_{k=1}^{\infty}\in H$
\begin{align*}
\E|v_{N}(\tau)-v_{M}(\tau)|^{2}_{H_0}&=\E\left|\sum_{k=1}^{\infty}\sum_{j=N+1}^{M}\eta_{j} \displaystyle\int_{0}^{\t}\langle{S}_A(t-s)\varphi_{j}(s),h_k\rangle_{H} h_k ds\right|^{2}\\
&=\sum_{k=1}^{\infty}\E\left|\sum_{j=N+1}^{M}\eta_{j}\displaystyle\int_{0}^{\t}\langle{S}_A(t-s)\varphi_{j}(s),h_k\rangle_{H}ds\right|^{2}\\
&\leq\sum_{k=1}^{\infty}\sum_{j=N+1}^{M}\displaystyle\E(\eta_{j})^2\int_{0}^{\t}\left|\langle{S}_A(t-s)\varphi_{j}(s),h_k\rangle_{H}\right|^{2}ds\\
&\leq\sum_{j=N+1}^{M}\displaystyle\E(\eta_{j})^2\int_{0}^{\t}\sum_{k=1}^{\infty}\left|\langle{S}_A(t-s)\varphi_{j}(s),h_k\rangle_{H}\right|^{2}ds\\
&\leq\E\sum_{j=N+1}^{M}\displaystyle\int_{0}^{\t}|S_A(t-s)\varphi_{j}(s)|^{2}_{H}ds\\
&\leq\E\sum_{j=N+1}^{M}\displaystyle\int_{0}^{\t}|\varphi_{j}(s)|^{2}ds\displaystyle\int_{0}^{\t}|S_A(s)|^{2}_{\calL({H})}ds
\end{align*}
Then $\{v_{N}\}_{N=1}^{\infty}$ is a Cauchy sequence in $L^{2}_{\calF_\tau}(\Omega,H_0)$. Hence the following quantity
\begin{equation*}
  \int_{0}^{\t}S_A(t-s)d{W}(s)=\lim_{N\to +\infty}\int_{0}^{\t}{S}_A(t-s)d\mathscr{W}_{N}(s)
\end{equation*}
is well defined in $L^{2}_{\calF_\tau}(\Omega,H_0)$. Moreover, by the martingale property of the stochastic integral (see \cite{Da-Za}), we have that
\begin{align*}
  \E\left(\sup_{0\leq t\leq\tau}\|v_{N}(t)-v(t)\|^{2}_{{H}_0}\right)\to 0\quad\text{as}\quad N\to +\infty,
\end{align*}
where $v$ is defined by
\begin{align*}
  v(t):={S}_A(t)\xi+\int_{0}^{t}{S}_A(t-s)d{W}(s),\quad \P-a.s.
\end{align*}
Consequently, one can consider a subsequence $v_{N_{k}}$ converging $\P$-a.s, and uniformly in $[0,\tau]$. Now, as the first derivative operator is closed and under assumption \eqref{Assumption1.2}-(b), the following sequence
\begin{align}\label{nec-condi}
\frac{d}{dt}\int_{0}^{t}{S}_A(t-s)d\mathscr{W}_{N}(s)
\end{align}
converges $\P$-a.s, and the following process
\begin{equation*}
  \int_{0}^{\t}\frac{d}{dt}S_A(t-s)d{W}(s)=\lim_{N\to +\infty}\int_{0}^{\t}\frac{d}{dt}{S}_A(t-s)d\mathscr{W}_{N}(s)
\end{equation*}
is well defined. By uniqueness of the solution we end up with
\begin{align*}
X(t)&=\lim_{N\to +\infty}v^{\prime}_{N}(t)\\
&={T}_{0}(t)\xi+ \int_{0}^{\t}\frac{d}{dt}S_A(t-s)d{W}(s)\\
&={T}_{0}(t)\xi+\lim_{N\to +\infty}\lim_{\lambda\to +\infty} \int_{0}^{t}{T}_{0}(t-s)\lambda R(\lambda,{A})d\mathscr{W}_{N}(s),
\end{align*}
for any $t\in [0,\tau]$ and $\P$-a.s.
\end{proof}
\section{Application to parabolic  equations with white noise boundary conditions}\label{sec:3}
In some applications the noise can effect the evolution of a system only through the boundary of a region. So, let $H$ and $U:=\partial H$ be two separable Hilbert spaces. Consider the following evolution equation with boundary white-noise
\begin{align}\label{Boundary-pb}
\begin{split}
\frac{\partial X}{\partial t}(t)&=A_{m}X(t), \quad t\in[0,\tau],\\
\gamma X(t)&=\dot{W}(t),\quad  t\in(0,\tau],\\
X(0)&=\xi.
\end{split}
\end{align}
Here $A_{m}:D(A_{m})\subset H\to H$ is an operator with partial derivatives, $\gamma:D(A_m)\subset H\to U$ is a linear operator defining boundary conditions, and $W(t)$, $t\geq0$, is a cylindrical Wiener process taking values in a space of distributions on $U$, defined as in Section \ref{sec:00}. Note that the solution of \eqref{Boundary-pb} is called to be an Ornstein-Uhlenbeck process. Now, as the previous section we fix  $\{\varphi_{j}\}_{j=1}^{\infty}$ an orthonormal and complete basis in $L^{2}(0,\tau;{U})$, consisting of smooth functions and consider approximation equation of \eqref{Boundary-pb}
\begin{align}\label{approx-pb2}
\begin{split}
\frac{\partial X_{N}}{\partial t}(t)&=A_{m}X_{N}(t), \quad t\in[0,\tau],\\
\gamma X_{N}(t)&=\dot{W}_{N}(t),\quad  t\in(0,\tau],\\
X_{N}(0)&=\xi.
\end{split}
\end{align}
We firstly reformulate the approximate boundary problem \eqref{approx-pb2} into stochastic Cauchy problem with non-dense domain. In order to do so, we introduce the product space
\begin{equation*}
  \mathscr{H}=U\times H,
\end{equation*}
which is a Banach space equipped with usual product norm. The following operator and process are important for the reformulation of the boundary white noise system \eqref{Boundary-pb} to a stochastic Cauchy problem with a non-dense domain
\begin{align*}
& \mathscr{A}\begin{pmatrix}0\\\phi\end{pmatrix}:=\begin{pmatrix}-\gamma\phi\\A_m\phi\end{pmatrix},\quad D(\mathscr{A}):=\left\{0_{{U}}\right\}\times D(A_{m}),\cr
& \mathscr{W}_{N}(t):=\begin{pmatrix}W_{N}(t)\\0\end{pmatrix},\quad t\in [0,\tau].
\end{align*}
By identifying $X_{N}(t)$ to $Z_{N}(t)=\begin{pmatrix}0\\X_{N}(t)\end{pmatrix}.$ The stochastic boundary problem \eqref{approx-pb2} can be rewritten as the following abstract stochastic Cauchy problem
\begin{align}\label{Trans-1}
\begin{cases}dZ_{N}(t)&=\mathscr{A} Z_{N}(t)dt+d\mathscr{W}_{N}(t),\quad t\in [0,\tau], \cr Z_{N}(0)&=(\begin{smallmatrix}0\\\xi\end{smallmatrix}).\end{cases}
\end{align}
One may observe that the operator $\mathscr{A}$ is turns out to be non-densely defined as $\overline{D(\mathscr{A})}\neq \mathscr{H}$. This reformulation allows us to make the following remark.
\begin{remark}
The stochastic approximative system \eqref{approx-pb2} is well posed if and only if the stochastic Cauchy problem \eqref{Trans-1} is well posed.
\end{remark}
To investigate the well posedness of stochastic Cauchy problem \eqref{Trans-1}, we construct the part of $\mathscr{A}$ in $\overline{D(\mathscr{A})}=\left\{0_{\mathcal{U}}\right\}\times H=:\mathscr{H}_{0}$, denoted by $\mathscr{A}_{0}:D(\mathscr{A}_{0})\subset\mathscr{H}_{0}\to\mathscr{H}_{0}$ which is a linear operator given by
\begin{align*}
 \mathscr{A}_{0}\begin{pmatrix}0\\\phi\end{pmatrix}=\begin{pmatrix}-\gamma\phi\\A_m\phi\end{pmatrix},\quad D(\mathscr{A}_{0}):=\left\{\left(\begin{smallmatrix}0\\\phi\end{smallmatrix}\right)\in \left\{0_{\mathcal{U}}\right\}\times D(A_{m})):\gamma\phi=0\right\}.
\end{align*}
We need the following assumptions in order to apply the main result (see Theorem \ref{main-result2}).
\begin{itemize}
  \item The operator $\mathscr{A}_{0}$ generates an analytic $\mathcal{C}_{0}$-semigroup on $\mathscr{H}_{0}$ denoted by\\ $(\mathscr{T}_{0}(t))_{t\geq0}$.
  \item There exist $\omega\in\mathbb{R}$ and $p^{*}\in(0,+\infty)$ such that $(\omega,+\infty)\subset\rho(\mathscr{A})$ and the resolvent set of $\mathscr{A}$ satisfies
      \begin{equation*}
        \lim_{\lambda\to +\infty}\lambda^{\frac{1}{p^{*}}}\|R(\lambda,\mathscr{A})\|_{\mathscr{H}}<+\infty.
      \end{equation*}
  As concluded $\mathscr{A}$ generates an integrated semigroup on $\mathscr{H}$ denoted by $(\mathscr{S}_A(t))_{t\geq0}$ (see section \ref{sec:1}).
 \begin{equation*}
  \int_{0}^{\tau} \|\mathscr{S}_{A}(r)\|^{2}_{\calL_{2}(\mathscr{H})}dr<\infty.
  \end{equation*}
 \item[(b)] As $t\to \mathscr{S}_A(t)$ is continuously differentiable, we assume that
 \begin{equation*}
   \int_{0}^{\tau} \left\|\frac{d\mathscr{S}_{A}(r)}{dr}\right\|^{2}_{\calL_{2}(\mathscr{H})}dr<\infty.
\end{equation*}
 $\calL_{2}(\mathscr{H})$ stands for the space of Hilbert-Schmidt operators acting from $\mathscr{H}$ into $\mathscr{H}$.
\end{itemize}
Under above assumptions and according to Lemma \ref{lemma-result1} we have for any $(\begin{smallmatrix}0\\\xi\end{smallmatrix})\in\mathscr{H}_{0}$, the mild solution of stochastic boundary noise problem \eqref{Trans-1} satisfies
\begin{align*}
Z_{N}(t)=\mathscr{T}_{0}(t)Z_{N}(0)+\lim_{\lambda\to +\infty} \int_{0}^{t}\mathscr{T}_{0}(t-s)\lambda R(\lambda,\mathscr{A})d\mathscr{W}_{N}(s),
\end{align*}
for any $[0,\tau]$ and $\P$-a.s.\\
By identifying $X(t)$ to $Z(t):=\begin{pmatrix}0\\X(t)\end{pmatrix}$ and $\mathscr{W}(t):=\begin{pmatrix}W(t)\\0\end{pmatrix}$, the SDE  \eqref{Boundary-pb} can be rewritten as the following stochastic Cauchy problem
\begin{align}\label{Trans-cau}
\begin{cases}dZ(t)&=\mathscr{A} Z(t)dt+d\mathscr{W}(t),\quad t\in [0,\tau], \cr Z(0)&=(\begin{smallmatrix}0\\\xi\end{smallmatrix}).\end{cases}
\end{align}
If the above assumptions hold, then by Theorem \ref{main-result2} we get that for any $(\begin{smallmatrix}0\\\xi\end{smallmatrix})\in\mathscr{H}_{0}$ the mild solution satisfies
\begin{align}\label{Mild-solution}
Z(t)=\mathscr{T}_{0}(t)(\begin{smallmatrix}0\\\xi\end{smallmatrix})+\lim_{\lambda\to +\infty} \int_{0}^{t}\mathscr{T}_{0}(t-s)\lambda R(\lambda,\mathscr{A})d\mathscr{W}(s),
\end{align}
for any $t\in [0,\tau]$ and $\P$-a.s.
\begin{remark}
These kind of result are already available for stochastic parabolic equations in extrapolation spaces (a larger space than $\mathscr{H}$), see e.g \cite[Theorem 2.1, Theorem 2.2]{Daprato-93}. Here we have offered a different proof for this particular case. However, in \eqref{Mild-solution} the solution takes their values in the state space $\mathscr{H}_0$.
\end{remark}
\begin{example}
Let $\mathscr{O}\subset\mathbb{R}^{n}$ be an open bounded smooth domain with boundary $\partial\mathscr{O}=:\Gamma$. Consider the following stochastic parabolic problem
\begin{align}\label{Neumann-Example}
\begin{cases}
\frac{\partial X}{\partial t}(t,x)&=\frac{\partial^{2}X}{\partial x^{2}}(t,x),\quad  t\in[0,\tau],x\in\mathscr{O}\\
\partial_\nu X(t,x)&=\dot{W}(t),\quad t\in(0,\tau],x\in\Gamma\\
X(0)&=\xi.
\end{cases}
\end{align}
Here $\nu=\nu(x)$ is the unit normal vector of $\Gamma$ pointing towards the exterior of $\mathscr{O}$, the initial condition $\xi\in
L^{2}(\Omega\times\mathscr{O})$. Note that the trace mapping $\gamma_{\Gamma}$ is continuous from the Hilbert space $H^{1}(\Gamma)$ into the Hilbert space $H^{\frac{1}{2}}(\Gamma)$, for more details on the measure of boundary we refer to \cite{Droniou}.\\
The boundary noise $W$ taking values in $L^{2}(0,\tau;H^{\frac{1}{2}}(\Gamma))$ and as pointed out in Section \ref{sec:1} the process $W$ can be represented as a formal series. In order to state the problem in our abstract setting, we introduce the product space
\begin{equation*}
\mathscr{H}:=H^{\frac{1}{2}}(\Gamma)\times L^{2}(\mathscr{O}),
\end{equation*}
as well as we select the following operator
\begin{align*}
D(\mathscr{A}):=\left\{{0}\right\}\times W^{2,2}(\mathscr{O}),\quad \mathscr{A}\begin{pmatrix}0\\\phi\end{pmatrix}:=\begin{pmatrix}-\partial_\nu\phi\\\Delta\phi\end{pmatrix}.
\end{align*}
One may observe that
\begin{equation*}
\overline{D(\mathscr{A})}=\left\{{0}\right\}\times L^{2}(\mathscr{O})=:\mathscr{H}_{0}.
\end{equation*}
With these notations the system \eqref{Neumann-Example} takes the form of the stochastic Cauchy problem \eqref{Cauchy}. Thus to show the well-posedness of system \eqref{Neumann-Example} it suffices to check the assumptions of Theorem \ref{main-result2}. According to the results obtained by Agranovich et al \cite{Agra-Denk-97}, the operator $\mathscr{A}$ satisfies the assumption \eqref{Assumption1.1}-(a). Thus it generates a nondegenerate exponentially bounded integrated semigroup denoted as $(\mathscr{S}_{\mathscr{A}}(t))_{t\geq0}$. The map $t\to \mathscr{S}_{\mathscr{A}}(t)$ is continuously differentiable on $(0,+\infty)$ and $(\frac{\partial \mathscr{S}_{\mathscr{A}}(t)}{\partial t})_{t\geq0}:= (\mathscr{T}_{0}(t))_{t\geq0}$ forms an analytic $\mathcal{C}_0$-semigroup on $\mathscr{H}_{0}$ which implies assumption \eqref{Assumption1.1}-(b) (see \cite[Proposition 3.4.3]{Magal-Ruan}). Now, according to \cite[Lemma 3.7]{Ducrot-Magal-Prevost} and \cite[Lemma 2.2.10]{Magal-Ruan} we have for each $x\in \mathscr{H}$ and $\beta\in(\frac{1}{4},\frac{1}{2})$
\begin{align*}
\frac{d}{dt}\mathscr{S}_{\mathscr{A}}(t)x&=(-\mathscr{A}_0)^{\beta}\mathscr{T}_{0}(t)(-\mathscr{A})^{-\beta}x\\
&=(\lambda_k)^{\beta}e^{-\lambda_k t}(\lambda_k)^{-\beta}x,
\end{align*}
where $(\lambda_k)_{k\geq1}$ denoted the eigenvalues of $-\mathscr{A}_0$. Then
\begin{align*}
\displaystyle\int_{0}^{\tau}\left\|\frac{d \mathscr{S}_{\mathscr{A}}(r)}{dr}\right\|_{\calL_2(\mathscr{H})}^{2}dr&=\displaystyle\int_{0}^{\tau} Tr ((-\mathscr{A}_0)^{2\beta}e^{2r A_0}(-\mathscr{A})^{-2\beta}) dr\\
&=\sum_{k=1}^{\infty}\frac{1-e^{-2\lambda_k\tau}}{2\lambda_k}.
\end{align*}
Since $\lambda_k\sim ck^{\frac{2}{n}}$, the above series convergent if and only if $n=1$. Thus assumption \eqref{Assumption1.2} holds. By Theorem \ref{main-result2}, the parabolic stochastic boundary problem \eqref{Neumann-Example} admits a unique integrated solution.
\end{example}
\begin{remark}\label{Dirichlet-booundary-condition}
\begin{itemize}
  \item [$(i)$] Note that if the operator $\mathscr{A}_{0}$ is sectorial and $\mathscr{A}$ is almost sectorial with $\omega_{\mathscr{A}}<0$. It has been shown in \cite[Lemma 3.8]{Ducrot-Magal-Prevost} that if $\gamma$ is a Neumann Boundary operator, then we have for any $p\in [1,4)$
\begin{equation*}
\displaystyle\int_{0}^{\tau}\left\|\frac{d \mathscr{S}_{\mathscr{A}}(r)}{dr}\right\|_{\calL(\mathscr{H})}^{p}dr <\infty.
\end{equation*}
\item [$(ii)$] Now, let us consider the problem \eqref{Boundary-pb} with inhomogeneous Dirichlet conditions. As above, we reformulate the boundary Dirichlet problem with white noise boundary data as an evolution equation in the following product space
\begin{equation*}
\mathscr{H}:=H^{\frac{3}{2}}(\Gamma)\times L^{2}(\mathscr{O}).
\end{equation*}
To prove the well-posedness of Cauchy Dirichlet problem, it suffices to check the assumptions of Theorem \ref{main-result2}. According to \cite[Theorem 2.1]{Agra-Denk-97}, the assumption \eqref{Assumption1.1} holds for $p^*=4$, i.e there exist $\omega\in\mathbb{R}$ such that $(\omega,+\infty)\subset\rho({\mathscr{A}})$ and the resolvent set of $\mathscr{A}$ satisfies for $p^*=4$
      \begin{equation*}
        \lim_{\lambda\to +\infty}\lambda^{\frac{1}{p^{*}}}\|R(\lambda,\mathscr{A})\|_{{H}}<+\infty.
      \end{equation*}
Which means that the associated integrated semigroup $(\mathscr{S}_{\mathscr{A}}(t))_{t\geq0}$ is nondegenerate exponentially bounded integrated semigroup (see \cite[Proposition 3.4.3]{Magal-Ruan}). It remains to verify the estimation \eqref{main-result2}. As shown in \cite[Lemma 3.8]{Ducrot-Magal-Prevost} the estimate in question holds if and only if
\begin{equation*}
\displaystyle\int_{0}^{\tau}\left\|\frac{d \mathscr{S}_{\mathscr{A}}}{dr}(r)\right\|_{\calL(\mathscr{H})}^{p}dr <\infty,
\end{equation*}
for any $p\in [1,\frac{4}{3})$.
\end{itemize}
\end{remark}

\end{document}